# Criteria for Solar Car Optimized Route Estimation

Mehrija Hasicic, Damir Bilic and Harun Siljak[1]

Electrical and Electronics Engineering Department, International Burch University,

Francuske revolucije bb, Sarajevo, Bosnia and Herzegovina.

Email: {name.surname}@ibu.edu.ba

**Abstract.** This paper gives a thorough overview of Solar Car Optimized Route Estimation (SCORE), novel route optimization scheme for solar vehicles based on solar irradiance and target distance. In order to conduct the optimization, both data collection and the optimization algorithm itself have to be performed using appropriate hardware. Here we give an insight to both stages, hardware and software used and present some results of the SCORE system together with certain improvements of its fusion and optimization criteria. Results and the limited applicability of SCORE are discussed together with an overview of future research plans and comparison with state-of-the-art solar vehicle optimization solutions.

---

1 Corresponding author. Email: harun.siljak@ibu.edu.ba

## 1. Introduction

Development of navigation systems has been an important topic in optimization once portable electronic devices were feasible [1, 2] and route selection and optimization have been the vital part of it, aiming at fuel consumption reduction [3] and driver satisfaction, which is in general a multi-objective optimization problem [4]. With the advent of electric and autonomous cars, attention in development of optimization algorithms turned to them, utilizing properties of these new vehicles [5-7]. Navigation for autonomous vehicles also allows use of algorithms previously developed for mobile robotics [8].

Geographic Information System (GIS) integrates various data types, many of which are instrumental in navigation, leading to extensive use of GIS in route optimization [9]. From our perspective, it is also important to note utilization of GIS (solar radiation maps) in solar energy management [10, 11] as it allows us to use GIS for our optimization as well as an input provider.

In the spirit of optimization with respect to fuel consumption, power management optimization in electric and solar cars has been investigated in theory and practice [12], and for solar hybrid cars the power management schemes focus on the question of switching energy sources and optimization of resources [13]. Sunshine forecast as an optimization input has been recently introduced [14, 15] and used mainly for parking planning.

This paper presents Solar Car Optimized Route Estimation (SCORE), a novel route optimization system based on proposing sunniest routes and sunniest parking spots, therefore utilizing the options of charging while driving and charging on parking lots. The data on solar irradiance for routes and parking spots is a fusion of previously collected, real time and forecasted data.

Describing the whole process from data collection to route selection, this paper provides both theoretical and practical treatise of SCORE, giving a general structure and the real world implementation of it. The paper itself is an extension of the work under same title presented in MECO 2016 conference [16], presenting the implementation challenges and solutions [17]. In addition to previous work, this paper also extends the analysis of cost functions used in SCORE for selection of routes and parking places, as well the mathematical model of data fusion, together with more details on the results.

## 2. State of the Art

Route planning and selection for road vehicles has been a subject of interest for decades, with the first commercial digital map navigation system appearing more than 30 years ago [18]. Since then, various route planning systems have been proposed, based on different algorithms and inputs.

Dijkstra's algorithm has been the simplest algorithm for implementation and serves as a benchmark for storage and time consumption in route planning for cars [19] when compared with other algorithms, from bidirectional and A* search to special goal-driven, hierarchical and bounded-hop algorithms [20]. Although Dijkstra's algorithm is the slowest option in big networks, it can still be used as a proof of concept.

While in the beginning the route planning and selection systems had a single objective: namely, fuel consumption or journey time minimization, soon after their inception multi-objective optimization models were developed, often using artificial intelligence techniques to combine different goals [21]. These goals often include personalized choices of drivers [22] and their personal attitude towards possible routes [23, 24].

Parking selection has been studied extensively as well [25]. It has been modeled as a multi-input problem measuring the utility of a parking space by accounting for availability, driving duration, walking distance to the destination, parking cost, traffic congestion, etc [26].

In terms of solar car optimization, the power management techniques mentioned in the introduction (namely [13]) have been extended to minimise total energy consumption by planning speed on parts of the path differently exposed to the sun [27]. This work, published at the same time as the first SCORE results [16] builds up on closely related work on solar powered robots [28].

In [29] authors propose a solar race car optimization based on weather forecast and velocity profile, determining the need for acceleration and deceleration throughout the race course in order to maximise the average velocity. Similar task is done in [30] as well, but the latter also includes somewhat more complex weather model and solar position algorithm. The weather model is a random walk around expected irradiance curve, while solar position is determined through the NREL (National Renewable Energy Laboratory) model.

In comparison to these solutions, SCORE has somewhat different purpose. Aiming at offering a framework for optimized city and intercity travel, SCORE's model relies on data which can be collected in a more regular fashion than the data for a race track [30]. Moreover, SCORE's output is the route, unlike the outputs of power management optimization systems having the velocity profile or engine-generator power trajectory (in case of hybrid cars) as the output [13, 29, 30]. This does mean that SCORE should be extended so it covers the velocity control and/or hybrid vehicle power switching, but at this point the focus is on route selection. This is along the lines of the idea in [27], where the path selection is followed by speed profile selection.

## 3. SCORE System Description

SCORE system consists of three separable parts, indicated in Fig. 1, namely:

1. Mobile sensor data transmitter, transmitting solar irradiance data through wireless channel in real time from the roads. Although we will use the term mobile sensors throughout this paper, they can be stationary as well, placed at selected places by the road. When mobile, these transmitters are not necessarily placed on solar cars using SCORE as their navigation system. They can be placed on fossil fuel and electric cars as well. Preferably, cars carrying the sensors would be often in motion, covering a large area (e.g. taxis, public transportation).
2. Server for data fusion, collecting readings transmitters send from the field and third party sources, processing them and combining with offline data (which can include weather forecast and historic readings, as the next section will show) and allowing the car computer clients to fetch the processed data in appropriate matrix format.
3. Embedded car computer client in the solar car, taking the processed data from server's cloud service and customize it on its own by using readings from its own sensors. Built-in light sensor can be used for normalization of data, and electric measurements from the car can be utilized for state estimation. Finally, the user can enter the destination and obtain the proposed route, which should dynamically change based on weather updates.

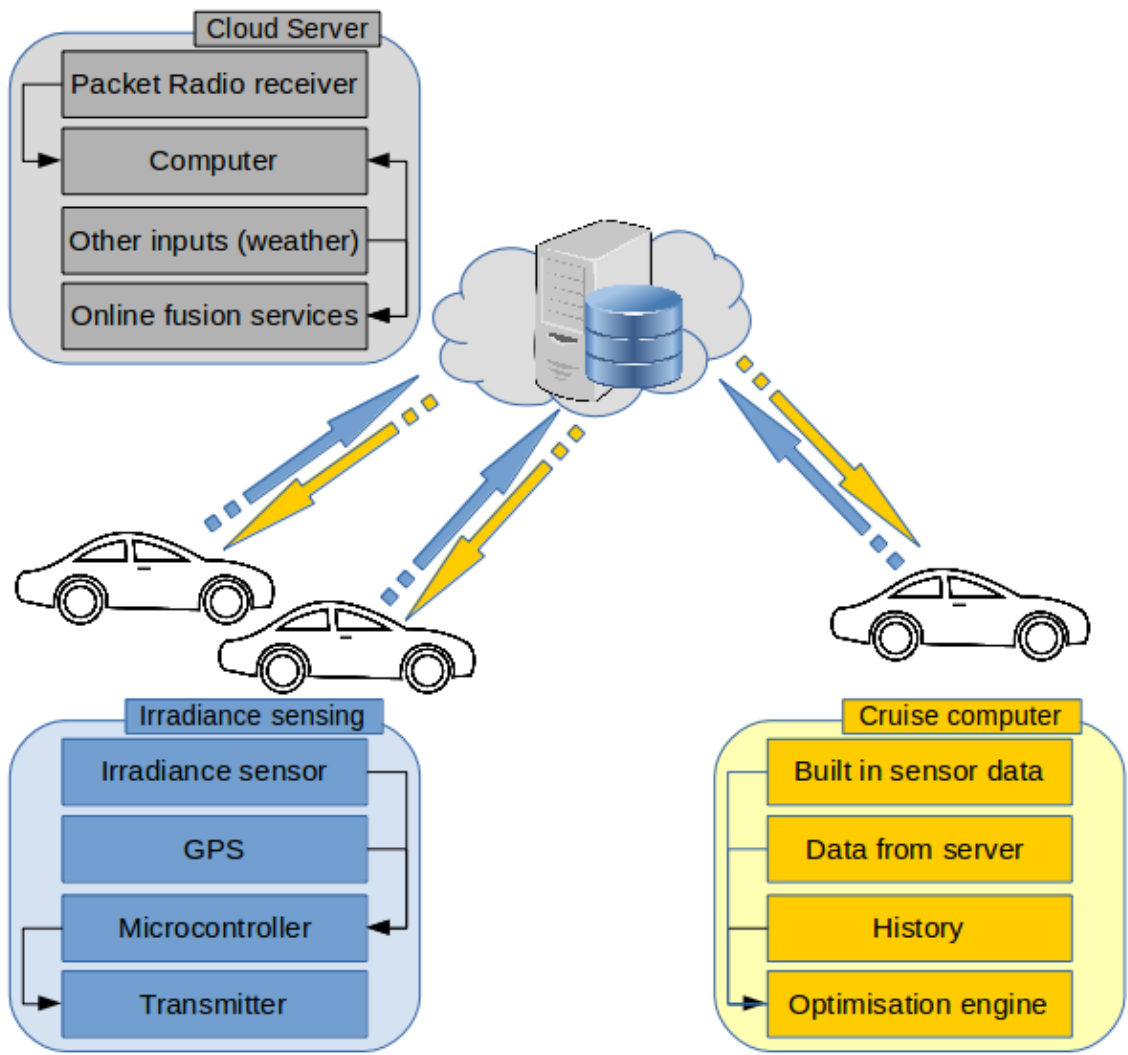

Figure 1: Overall structure of SCORE

## 4. Data Collection

### 4.1. Theoretical Considerations of Irradiation Data Analysis

In order to select routes with highest solar energy gain, SCORE system has to have relevant solar irradiation data. Since it is not possible to always have up to date data in real time for every road segment considered by the algorithm, it is important to use different sources of information. In this work, we have divided the irradiation data into two categories:

1. Online data, gathered by the mobile sensor data transmitters and updated in regular fashion. The data for each location is represented by real numbers between 0 (no irradiance) and 1 (maximum irradiance) with a timestamp for data sample collection. In this paper, timestamps are integers denoting hours starting from a reference time (beginning of the year).
2. Offline data, generated using numerical sunshine forecast, CAD (computer aided design) and GIS models for prediction of solar irradiance for a particular location. CAD data is generated from CAD street models and simulating sun movement [31, 32], while GIS data is taken from the GIS services doing solar irradiance measurements for areas of interest [10]. This data is provided in an aggregated form by Google as well through their Project Sunroof for housing solar panel planning.

These two numerical values, denoted $r_{on}$ (normalized value of online data irradiance) and $r_{off}$ (normalized value of irradiance inferred from offline data) are combined for each geographical location (in the optimization part, we will refer to these as graph nodes) on the server. Details of this fusion will be discussed in a separate section.

### 4.2. Implementation of Sensor Data Collection and the Server

Mobile device developed within this project is compact and autonomous, which enables its placement on a vehicle moving through the city to collect irradiation data without customization of the car itself or its routes [17]. Once again it is emphasized that the vehicles carrying these mobile devices do not have to be vehicles using SCORE for navigation, i.e. traditional fossil fuel cars may serve the purpose of data collecting "crawlers".

While any wireless protocol could be used for transmission from these mobile devices, we propose the use of packet radio. Its easiest implementation is APRS (Automatic Packet Reporting System) which has been used to deliver GPS and sensor data to the terminal node, connected to the server. APRS has been used before for such purposes as well, in monitoring systems [33]. For the transmission, we have used amateur radio bands. Of course, different implementations of SCORE could use other frequency bands and/or proprietary protocols if needed.

One may argue that the data from the devices could be simply kept in memory and read at the end of the day, which would be appropriate if we were only interested in historic data and statistics. With real-time updates from the mobile devices, our servers always have fresh data and improve the relevance of SCORE's path selection for the clients, as the optimization algorithm relies more on the new data, using old data merely as a reference.

Large computing power is not a requirement for the server (in our case, it was Raspberry Pi 2) and many operations, depending on the implementation of SCORE, may be run on cloud as well. Storing fusion tables (graph matrices, as we will refer to them in the optimization part) on cloud enables both the clients and third parties with their smartphones, computers and dedicated hardware to access the processed data for their own uses.

Prototype of the server in our case receives the radio packed data, converts it from audio to text using a standard sound card, custom interface developed in [34] controlled by dedicated software (e.g. AGWMonitor). This data is merged with CAD and GIS data as suggested in previous section and placed on the cloud as an optimization input for clients (as described in the next section).

## 5. Optimization

### 5.1. Theoretical Considerations of Optimization

Dijkstra's algorithm [35] was a straightforward choice for route selection. The input to the algorithm is a graph (in computing terms, represented as a matrix) with positive weights of edges to be defined, and a set of nodes. In our case, the nodes are major crossroads. In a classical fuel optimization problem for fossil fuel cars, the weights of branches would be road lengths. However, in our case we account for the converted solar energy on the road segment as well. This value is calculated using the known parameters of our vehicle:

- 11 kW motor power
- 2x0.726 $m^2$ panel area
- 18% panel efficiency
- 957 $W/m^2$ received power per square meter with full (unity) irradiation (i.e. reflection of 30%)

Using these parameters, we get a new length of road segment which is shorter physically, as solar conversion compensates for some of the energy spent. These values, the new lengths, constitute the graph matrix forwarded to Dijkstra's algorithm.

In this work, we have used a first order approximation, taking the solar irradiation of the segments to be the arithmetic mean of irradiation in nodes. While it is a crude estimate, it is still satisfying the needs of prototype testing. Once a massive network of mobile sensors and transmitters is deployed, the need for approximation disappears.

As far as parking space selection is concerned, a straightforward algorithm was originally used, defining a cost function to be the quotient of irradiation divided by the distance of parking space from target destination. Taking the numerator or the denominator to non-unit powers can tweak the importance of distance or the irradiation in parking space selection, according to user's needs. An extension of this consideration is given in the next section.

These quotients and graph weights from previous discussion are calculated on the server and provided to the client through cloud, in order to reduce computation burden on the client architecture.

### 5.2. Implementation of the Optimization Client

While the prototype device for the solar car has been developed on a microcontroller, another testbed for algorithm tests was developed on a PC using MATLAB (MATLAB system works as an integrated server/client environment).

The car computer prototype is built on TI's ARM Cortex-M4F based TM4C123G LaunchPad because this development board had significantly more working memory than Arduino, our original prototyping platform. More memory was needed to keep the whole matrix representing the graph in the working memory of our embedded processor. We have noted that the matrix in question is sparse, so the implementation would benefit from optimized storage. However, having a relatively small (50x50) matrix in our case meant that the system will not encounter problems with dealing with the matrix in its raw format. The embedded system has a keyboard (used to enter the destination node) and display (showing the selected route) placed in the custom-built solar car. In this prototype, client receives the graph matrix through wireless or USB debug cable and performs the route selection (i.e. Dijkstra's algorithm) to provide the route (sequence of nodes) to the user, based on the desired destination node.

The full-fledged client will be using:

1. Sensor fusion data from the cloud
2. Measurements from the solar panels and battery
3. Built-in irradiation sensor measurements
4. Routes from user's history

In this prototype implementation, only the first set of inputs is used for simplicity, because it is the main source of data for SCORE. Measurement (2) is needed for feedback, measurements (3) are a corrective input (for normalization of data fetched from the cloud) and routes (4) are used to improve user experience.

## 6. Optimization Criteria and Data Fusion

Online and offline data fusion, previously mentioned in subsection 4.1 is a vital part of SCORE. It enables clients to have relevant data even if the field records are not up to date. Previously, the following formula was used for fused irradiance *r*, combining online data $r_{on}$ and offline data $r_{off}$:

$$r = a \cdot r_{on} + (1-a) \cdot r_{off}, \text{ where } a = \exp\left(\frac{-(t_{curr} - t_{meas})^2}{100000}\right) \tag{1}$$

Here we use timestamps mentioned before, determining how many hours have passed between the last online measurement for the position in question ($t_{meas}$) and current moment ($t_{curr}$). The denominator in *a* has been empirically chosen.

However, this model can be vastly improved. For instance, the following value of *a*

$$a = \exp(-m^2 - d) \text{ where } m = mod(t_{curr} - t_{meas}, 24), \ d = \lfloor \frac{t_{curr} - t_{meas}}{24} \rfloor \tag{2}$$

would take into account the time of the day (i.e. give preference to data collected at the same time of the day the planned trip takes place in) but also number of days since the data was collected. Here the mod (division remainder) is taken to be from the interval [-11,12] and not [0,23].

Three days of hourly solar irradiation data measurements taken from NREL [36] were used for testing of these online-offline data fusion models in the following experiment. Value of *r* at 4 PM of day 3 was supposed to be estimated based on the rest of the data collected before 4 PM (we only took into account the 31 hours with non-zero irradiation, i.e. day hours, 12 hours of night for each 24-hour cycle were left out).

For each of the 31 non-zero irradiation hours we took the irradiation value and timestamp difference from the 3<sup>rd</sup> day, 4 PM mark and computed *r* using (1) and (2), assuming that it is the last piece of data for particular location that has been collected. The $r_{off}$ value was taken to be a random variable from (0,1) interval, uniformly distributed. The trials were run 100 times, averaged and results are shown in Fig. 2.

The empirical choice of parameter in (1) to be 100,000 ensures data from a different season (several months ago) is not used for the fusion. However, it filters out the forecast completely in a few hours or a few days estimation window, as seen in the figure (dashed line) and relies solely on the historical data. This is why we substituted it with a smaller value (10) and ran the calculation (dot dash line). Convergence to the forecast value is significantly faster. Intermediate values between 100,000 and 10 were also tested and the results were equivalent to those seen for 100,000 as the influence of forecast was still too small on two day window width, and the parameter value of 10 completely discards data older than 6-8 hours and turns solely to forecast in absence of fresh data.

Now, let us examine the proposed fusion with mod 24 (represented by (2) and solid line in the graph). It converges to the forecast value fast, but diverges if the data is taken at the same time on the previous day, improving the quality of forecast. This effect also decays fast, so the measurements from the same time two days before are not taken as very significant. This analysis suggests that in case of unreliable forecast, model (2) can bring some improvement to the estimation of the route conditions. If the forecast is reliable, then there is no need for enhancement of this sort. Reliability of forecast should be subject to statistical analysis of the forecast vs. real state in a longer time period.

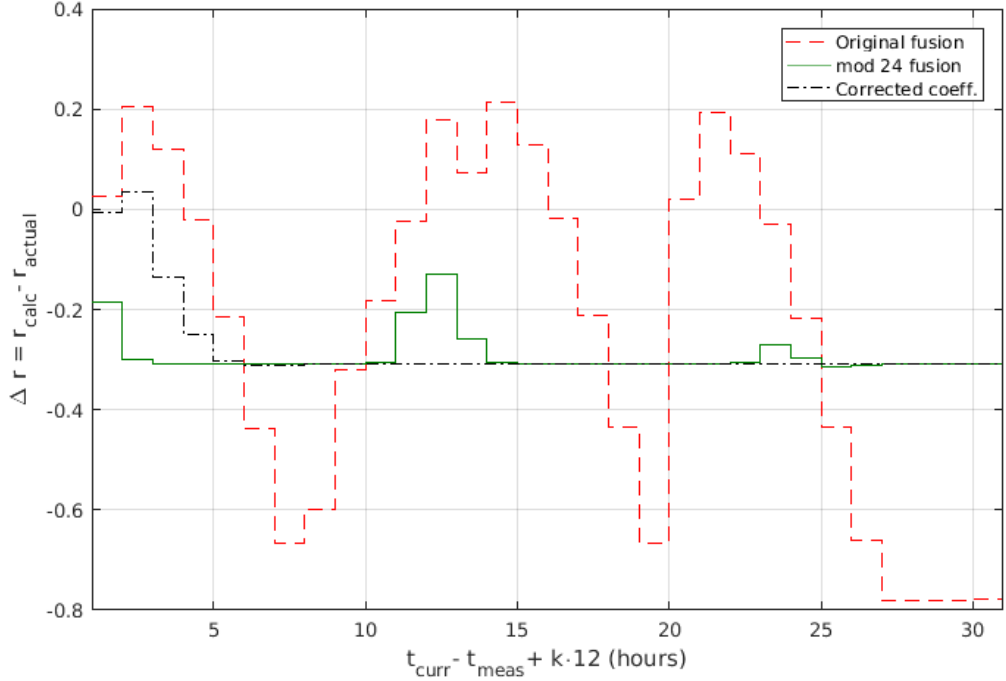

Figure 2: Error of three online-offline data fusion models predicting irradiation level with respect to the time of last data collection

On the other hand, in the optimization process the selection of parking spot has been left vaguely defined as well. While selection of route focuses on length (reduced length, to be precise), selection of parking spaces can be a multicriterial optimization problem, as shown in the state of the art considerations. In our particular simplified case, we have proposed earlier a quotient of irradiance and walking distance to be the optimization criterion [16] (as recalled in subsection 5.1).

However, simply posing it as *r/d* does not work, as the same percentage change in distance and irradiance should not have the same effect on place selection. Furthermore, as noted in [37], people will walk 100 to 500 m from their parking spot if it is necessary, not more. Hence we propose a new criterion of the form

$$f = \frac{r}{d^b + k}$$

where *r* and *d* are irradiance and walking distance from the spot to the target point, *b* is an exponent taken to be 1 in our analysis and *k* is a length parameter (in our analysis taken to be 200 m). Parameters *b* and *k* can be tweaked according to user's preference and prior experience. For a more detailed analysis, we suggest the criterion to be expanded taking into account other relevant factors listed in section 2.

## 7. Results

Application of SCORE in different scenarios leads to a conclusion that solar conversion has a rather low influence in route selection. Solar conversion results in 2-3% reimbursement of energy spent under the maximum (unity) irradiation condition, which means that the effect of optimization cannot be seen in city environment.

While the numerical simulation results show no effect of solar conversion shares of 10-15% in medium sized urban areas, this does not mean there would not be any benefit of SCORE in that case: it just suggests that a higher share would significantly change the route selection and show all advantages of SCORE.

In case of the solar conversion effect increase to 25% for maximal irradiation, route selection could be affected by it even in city conditions. This of course asks for a significant change in solar car features listed earlier. An example of two such routes is shown in Figure 3: path from point A (bottom left) to point B (top right) is changed from black in case of low solar conversion effect to red in case of high conversion effect (25%). Green part is common for both scenarios.

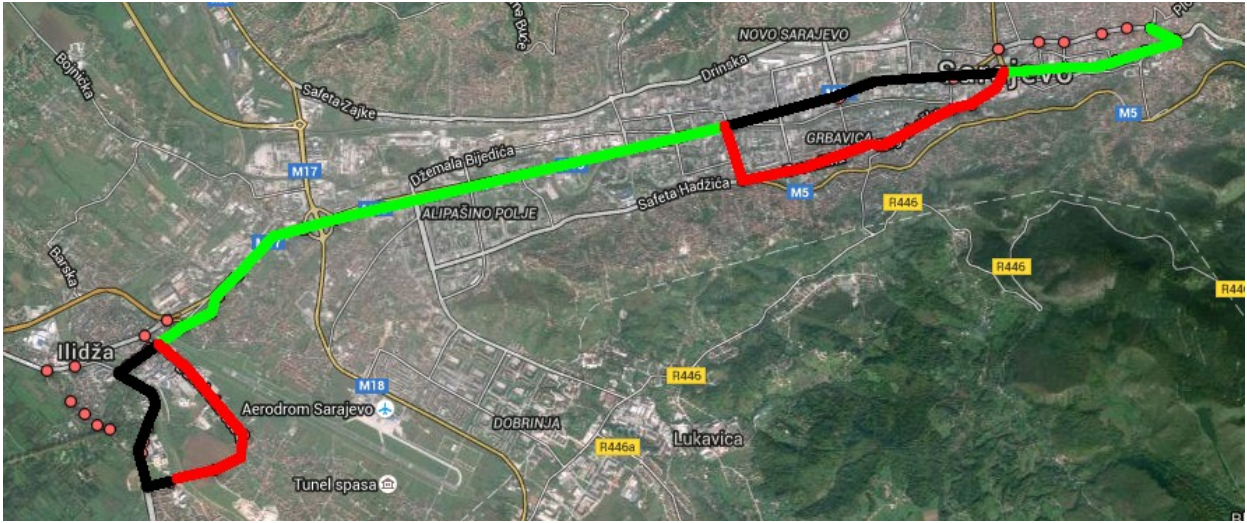

Figure 3: Influence of solar panel conversion rates on longer path choice

Parking selection gives the best results, even with the low conversion ratio, suggesting that a sunny parking spot has to have a priority over a close parking spot, and that the corresponding cost function should be tweaked in this sense (high power in the denominator which would lower the influence of distance). These conclusions are drawn from simulations showing the effect of long-term stay at a single location and the cost function being the ratio of irradiation and distance of the parking space from the target point. However, this cost function could be bilinear as well, adding offset constants in numerator and denominator, which in turn can be used to tune the algorithm and provide better user experience.

8. Discussion

Results presented have shown limited usability of SCORE system with the current solar car design. However, we suggest the need for change of car features, namely the solar panels. In particular, efficiency increase would greatly influence applicability of SCORE. Under the assumption of the impossible 100% efficiency and coverage of the whole top side of the car with panels, we may achieve necessary effects in city environment. While the current state of technology cannot offer high efficiency [38], improvement in manufacture [39] and reduction of losses [40] will help in increase of SCORE's applicability.

The results presented come from a real solar car and numerical simulations. However, with the modern hybrid cars with multiple sources of power, the question of power management and SCORE is left unanswered. Also we need to emphasize that effects of congestions, traffic jams and speed changes haven't been accounted for in this work. Their effect would be positive on the applicability of SCORE, increasing the converted solar energy ratio.

## 9. Conclusions

As a simple optimization of solar/hybrid car routes based on energy saving, SCORE may be an applicable solution for both solar cars and mobile robots with solar panels. Major limitation of it, the low influence of solar conversion on overall energy balance may be mitigated with the advancement of panel production technology, addition of regenerative braking or change of the car design.

Future work will aim at development of a network of mobile and stationary data collectors, which can collect even more data, relevant for fossil fuel and electric cars as well, such as road quality, traffic congestion, pollution. A large network of these collectors, placed on frequent drivers such as delivery trucks, taxis, public transport, helps in providing relevant data on a daily basis, mitigating the lack of 3D models of streets and GIS data for some areas.

Design upgrades may include developing more complex models for online-offline data combination, possibly by using artificial intelligence tools, as well as adding new inputs to the optimization process and experimenting with other optimization algorithms instead of Dijkstra's. This in particular means addressing the problem of optimization as a dynamical programming problem, enabling route changes on the go. Incorporation of sensory voltage and current readings from the car will also enable the feedback on energy spent and generated. Dijkstra's algorithm is slow compared to other algorithms used today [20], so the alternatives will be pursued. This is the only shortcoming of the proposed SCORE system with respect to existing systems, as it is the first one to offer this criterion of optimization for solar vehicles.

Question of incorporation of SCORE with existing power management schemes in hybrid solar cars remains an interesting topic for future research as well.

Finally, it is worth noting that the use of APRS for communication in SCORE allows ham radio operators to participate in both data collection and utilization, which in turn can create a crowdsourcing environment which has been successfully used already [4].

## Acknowledgement

Authors wish to thank Sabahudin Vrtagic, Semir Sakanovic (International Burch University Sarajevo) and their team for providing the solar car used in this paper. Authors also wish to thank Professor Samim Konjicija (University of Sarajevo) and Ms. Minja Miladinovic (Jozef Stefan Institute, Ljubljana, Slovenia) for their help in building the prototypes of devices used in a priori data collection.